\documentclass[12pt]{article}
\usepackage{subfigure}
\usepackage{epsfig}
\usepackage{amsmath, amsthm, amssymb}

\newcommand{\be}{\begin{eqnarray}}
\newcommand{\ee}{\end{eqnarray}}

\theoremstyle{plain} 
\newtheorem{dummy}{Dummy}
\newtheorem{theorem}[dummy]{Theorem}

\theoremstyle{definition}

\theoremstyle{plain}

\newtheorem{lemma}[dummy]{Lemma}
\newtheorem{proposition}[dummy]{Proposition}

\newtheorem{claim}[dummy]{Claim}

\begin{document}

\date{}
\title{The Poisson equation on complete manifolds with positive spectrum and
applications}
\author{Ovidiu Munteanu \thanks{omuntean@math.columbia.edu, Columbia University, New York, NY 10027.}
\and Natasa Sesum \thanks{natasas@math.columbia.edu, Columbia University/University of Pennsylvania, New York, NY 10027;
partially supported by the NSF grant DMS-06-04657.}}
\date{}
\maketitle

\begin{abstract}
In this paper we investigate the existence of a solution to the Poisson
equation on complete manifolds with positive spectrum and Ricci curvature
bounded from below. We show that if a function $f$ has decay $f=O\left(
r^{-1-\varepsilon }\right) $ for some $\varepsilon >0,$ where $r$ is the
distance function to a fixed point, then the Poisson equation $\Delta u=f$
has a solution $u$ with at most exponential growth.

We apply this result on the Poisson equation to study the existence of
harmonic maps between complete manifolds and also existence of Hermitian
Einstein metrics on holomorphic vector bundles over complete manifolds, thus
extending some results of Li-Tam and Ni.

Assuming that the manifold is simply connected and of Ricci
curvature between two negative constants, we can prove that in fact the
Poisson equation has a bounded solution and we apply this result to the
Ricci flow on complete surfaces.
\end{abstract}

\section{Introduction}

We consider a complete manifold $M$ with a Riemannian metric $g$ on it. We
assume that $M$ has a positive spectrum, which means that $\lambda _1(M)>0$,
where $\lambda _1(M)$ denotes the greatest lower bound of the $L^2$ spectrum
of the Laplace operator on $M$. The second condition is that the Ricci
curvature is bounded from below by a negative constant.

In the first part of the paper we study the Poisson equation on $(M,g)$, 
\begin{equation}
\Delta u=f,  \label{eq-pois}
\end{equation}
where $f$ is a function on $M$ with a decay 
\[
|f(x)|=O(d(p,x)^{-1-\varepsilon }), 
\]
with $d(p,x)$ being the distance function on $M$ with respect to a fixed
point $p\in M$ and $\varepsilon >0$.

The Poisson equation on complete manifolds has been an extensive subject of
study of many people over the years. The case of complete manifolds with
nonnegative Ricci curvature has been more approachable, and in fact Ni, Shi
and Tam gave in \cite{N-S-T} necessary and sufficient conditions for
existence of solutions with certain growth rates. Of great importance for
their study were the estimates of Green's function proved by Li and Yau \cite
{L-Y}. In general though, when the assumption that Ricci curvature is
nonnegative is removed, one does not have such estimates and the problem
becomes considerably more difficult. Results that we know assume either some
type of integrability for the function $f$ in the right hand side of (\ref
{eq-pois}), or that the manifold $M$ has a certain structure at infinity.
For example, it is known that a solution of the Poisson equation exists if
the manifold $M$ has positive spectrum and the right hand side of the
equation is in $L^p\left( M\right) $ for some $1<p<\infty $, see \cite{Ni1}, 
\cite{S}. Interesting results about the Poisson equation can also be found
in the context of geometric scattering theory, see e.g. \cite{Ma}, \cite{M-V}%
. However, this only applies when the manifold has controlled asymptotic
geometry.

Recall that our setting in this paper is that $M$ has positive spectrum and
Ricci curvature bounded below. It is well known that manifolds with positive
spectrum have exponential volume growth, see \cite{B}, \cite{L-W1},
therefore the assumption in \cite{Ni1} that $f\in L^p\left( M\right) $ with $%
p$ finite is quite restrictive. Our main result is that we can replace the
integrability condition with a mild pointwise polynomial decay and solve the
Poisson equation with a solution of at most exponential growth.

\begin{theorem}
\label{thm-poisson} Let $M$ be a complete Riemannian manifold of dimension $%
n.$ Assume that $M$ has positive spectrum and Ricci curvature bounded below
by a negative constant. Consider any bounded function $f$ having decay 
\begin{equation}
\left| f\left( x\right) \right| \leq \frac C{\left( 1+r\left( x\right)
\right) ^{1+\varepsilon }}.  \label{eq-decay}
\end{equation}
where $r\left( x\right) =d\left( p,x\right) $ and $\varepsilon >0.$ Then
there exists a solution $u$ of 
\[
\Delta u=f.
\]
Moreover, with respect to $p$ this solution has at most exponential growth,
i.e. there exist constants $A>0$ and $B>0$ such that for any $x\in M$ 
\[
\left| u\left( x\right) \right| \leq Ae^{Br\left( x\right) }.
\]
\end{theorem}

We will give some applications of Theorem \ref{thm-poisson}. One of them is
the existence of harmonic maps which are homotopic to a given map and in
that light we prove the following theorem.

\begin{theorem}
\label{thm-harmonic} Assume that $M$ is a complete Riemannian manifold with
positive spectrum and Ricci curvature bounded below. Let $N$ be a complete
manifold with nonpositive sectional curvature and $h$ a smooth map from $M$
to $N.$ If with respect to a fixed point $p\in M$ the tension field has
decay at infinity 
\[
\left| \left| \sigma \left( h\right) \right| \right| \left( x\right) \leq
\frac C{\left( 1+r\left( x\right) \right) ^{1+\varepsilon }},
\]
then there exists a harmonic map $u:M\rightarrow N$ such that $u$ is
homotopic to $h.$ Moreover, the homotopic distance between $u$ and $h$ has
at most exponential growth on $M.$
\end{theorem}

The second application is the existence of a Hermitian-Einstein metric on a holomorphic
vector bundle over a complete manifold.

\begin{theorem}
\label{thm-he} Let $M$ be a complete K\"{a}hler manifold with positive
spectrum and Ricci curvature bounded below. Assume that $E$ is a holomorphic
vector bundle of rank $k$ over $M$ with metric $H_0$ such that 
\[
\left| \left| \Lambda F_{H_0}-\lambda I\right| \right| \left( x\right) \leq
\frac C{\left( 1+r\left( x\right) \right) ^{1+\varepsilon }}.
\]
Then there exists a metric $H$ on $E$ such that $\Lambda F_H-\lambda I=0.$
\end{theorem}

Problems stated in Theorem \ref{thm-harmonic} and Theorem \ref{thm-he} have
been investigated in detail in \cite{Ni, Ni1, Ni-Ren} and we use the
arguments from those papers. The main difference is contained in the fact
that we are able to solve the Poisson equation under a lot milder
assumptions and therefore it allows us to prove quite general existence
theorems compared to the cited papers.

We also give an application to the Ricci flow on complete surfaces. We
consider a complete Riemann surface $M$ with the scalar curvature having a
decay 
\begin{equation}
|R(x)+1|=O(r(x)^{-1-\varepsilon }),  \label{eq-R-decay}
\end{equation}
where $\varepsilon >0$ and $r(x)=d(p,x)$, with $p\in M$ being a fixed point.
This means in some sense our manifold is close to a hyperbolic surface at
infinity. We will evolve such a metric by the Ricci flow equation, 
\begin{equation}
\frac \partial {\partial t}g=-(R+1)g,  \label{eq-rf}
\end{equation}
and study its long time existence and convergence as $t\to \infty $.

\begin{theorem}
\label{thm-rf1} Let $(M,g_0)$ be simply connected Riemann surface with its Ricci
curvature satisfying 
\[
-a^2\le \mathrm{Ric}_M\le -b^2<0,
\]
and (\ref{eq-R-decay}). Then the Ricci flow (\ref{eq-rf}) starting at $g_0$
exists forever and converges, as $t\to \infty ,$ to a complete metric of
negative constant curvature.
\end{theorem}

The organization of the paper is as follows. In section \ref{sect-prel} we
will give some preliminaries for the proofs of our main results. In section 
\ref{sect-poisson} we will prove Theorem \ref{thm-poisson}. The proof
includes various estimates on the growth of the Green's function. In section 
\ref{sect-appl} we will show how one can use Theorem \ref{thm-poisson} to
prove Theorems \ref{thm-harmonic}, \ref{thm-he} and \ref{thm-rf1}.

\section{Preliminaries}

\label{sect-prel}

In this section we will give some preliminaries about the Green's function
for a complete manifold $M$ with positive spectrum, harmonic maps between
two complete manifolds and Hermitian-Einstein metric on holomorphic vector
bundles over complete manifolds.

\bigskip

\textbf{Green's function:}

It is a classical result (\cite{S-Y}) that if the manifold has positive
spectrum then it is nonparabolic, that is, there exists a positive symmetric
Green's function $G$ on $M.$ Moreover, we can always take $G\left(
x,y\right) $ to be the minimal Green's function, constructed using
exhaustion of compact domains.

For any fixed $x\in M$ and any $0\leq \alpha <\beta \leq \infty $ denote by 
\begin{eqnarray}  \label{eq-level}
L\left( \alpha ,\beta \right) &=&\left\{ y\in M:\alpha <G\left( x,y\right)
<\beta \right\} .  \nonumber \\
l\left( s\right) &=&\left\{ y\in M:G\left( x,y\right) =s\right\},
\end{eqnarray}
the level sets of $G(x,\cdot)$. The dependence of $L\left( \alpha ,\beta
\right)$ and $l\left( s\right)$ on $x$ will always be assumed implicit in
the notation. We recall the following well known results about the level
sets of $G:$

\begin{enumerate}
\item[(i)]  Since $G$ is harmonic and with finite Dirichlet integral, 
\[
\int_{l\left( s\right) }\left| \nabla G\right| \left( x,y\right) dy
\]
is a finite number independent on $s.$

\item[(ii)]  $G$ has the following integral decay at infinity 
\begin{equation}
\int_{B_x\left( R+1\right) \backslash B_x\left( R\right) }G^2\left(
x,y\right) dy\leq \bar{C}e^{-2\sqrt{\lambda _1\left( M\right) }R},
\label{eq-green-decay}
\end{equation}

for any $R>2,$ and where $\bar{C}$ depends on $x$.

\item[(iii)]  The co-area formula and (i) give for any $\delta >0$ and $%
\varepsilon >0$ that 
\begin{gather}
\int_{L\left( \delta \varepsilon ,\varepsilon \right) }G^{-1}\left(
x,y\right) \left| \nabla G\right| ^2\left( x,y\right) dy=\int_{\delta
\varepsilon }^\varepsilon \frac 1tdt\int_{l\left( t_0\right) }\left| \nabla
G\right| \left( x,y\right) dy  \nonumber \\
=\left( \int_{l\left( t_0\right) }\left| \nabla G\right| \left( x,y\right)
dy\right) \left( -\log \delta \right) 
\end{gather}
\end{enumerate}

For the proof of these results see \cite{L-W}.

\bigskip
\bigskip

\textbf{Harmonic maps}

For a map between two manifolds $u:M\rightarrow N$ define the energy density
by 
\[
e\left( u\right) =g^{ij}\left( x\right) h_{\alpha \beta }\left( u\left(
x\right) \right) \frac{\partial u^\alpha }{\partial x^i}\frac{\partial
u^\beta }{\partial x^j},
\]
where $ds_M^2=g_{ij}dx^idx^j$ and $ds_N^2=h_{\alpha \beta }dy^\alpha
dy^\beta .$ Consider also $\sigma \left( u\right) $ the tension field of a
map $u$ i.e. 
\[
\sigma \left( u\right) =\left( \Delta u^\alpha +\Gamma _{\beta \delta
}^\alpha \frac{\partial u^\beta }{\partial x^i}\frac{\partial u^\delta }{%
\partial x^j}g^{ij}\right) \frac \partial {\partial y^\alpha }.
\]
We will say the map $u$ is harmonic if $\sigma \left( u\right) =0$.

\bigskip \bigskip

\textbf{Hermitian Einstein metric}

Let us recall some notation. Let $\left( E,H\right) $ be a holomorphic
vector bundle with Hermitian metric $H.$ Define the operator $\Lambda $ as
the contraction with the metric $g_{ij}$ on $M,$ i.e. 
\[
\Lambda \left( a_{i\bar{j}}dz^i\wedge d\bar{z}^j\right) =g^{i\bar{j}}a_{i%
\bar{j}}, 
\]
for any $\left( 1,1\right) $ form $a.$ The metric $H$ is called Hermitian
Einstein if the corresponding curvature satisfies 
\[
\Lambda F_H=\lambda I, 
\]
where $\lambda $ is some constant.

\section{The Poisson equation}

\label{sect-poisson}

The goal of this section is to prove Theorem 1. This result and the
estimates we prove here may have different applications for complete
manifolds.

Let $G(x,y)$ be the Green's function as in section \ref{sect-prel}. If we
manage to show that for every $x\in M$, 
\[
\left | \int_M G(x,y)f(y)\, dy \right | < \infty, 
\]
then a function $u(x)$ defined by the above integral solves (\ref{eq-pois}).
Moreover, we will show $u(x)$ has the right exponential growth as stated in
Theorem \ref{thm-poisson}. In order to control the above integral we need to
control the Green's function on $M$. That is the reason why we need to
establish certain growth estimates for $G(x,\cdot)$ first. Let us outline
the main steps in the proof of Theorem \ref{thm-poisson} and then give more
details below.

\begin{enumerate}
\item[(i)]  First we will show there are positive constants $A$ and $B$,
independent of $x$ so that for all $y\in \partial B_x(1)$, 
\[
A^{-1}e^{-Br(x)}<G(x,y)<Ae^{Br(x)},
\]
and that there exists positive constants $C$ and $b$ such that for any $x$ 
\[
\int_{l(s)}|\nabla G|(x,y)\,dy\leq Ce^{br(x)}.
\]

\item[(ii)]  It turns out that for any $x\in M$ we have $L(Ae^{Br(x)},\infty
)\subset B_x(1)$. Using (i) we show that 
\[
\int_{L(Ae^{Br(x)},\infty )}G(x,y)\,dy\le Ce^{Br(x)}.
\]
This estimate can be found in \cite{Yin} but we will include the proof of
it here for completeness.

\item[(iii)]  We will split the $\int_MG(x,y)f(y)\,dy$ into two parts, 
\[
\int_{L(Ae^{Br(x)},\infty )}G(x,y)f(y)\,dy\,\,\,\mbox{and}%
\,\,\,\int_{L(0,Ae^{Br(x)})}G(x,y)f(y)\,dy.
\]
To control the first one we will use the estimate in (ii) and to control the
second one we will use all the previous estimates and the decay for $f(x)$
on $M$. The proof for existence of $u$ is based on a level set argument. We
should point out that (i) and (ii) are used here only to prove the
exponential growth of the $u$, the solution of (\ref{eq-pois}).
\end{enumerate}

Let us first prove some results about the behavior of Green's function at
infinity. These results are of independent interest. The constants we obtain
will in general depend on the dimension of $M$, Ricci curvature lower bound
and the bottom of spectrum.

\begin{lemma}
\label{lem-green-point} There exist positive constants $A$ and $B$ such that
for any $x\in M$ and $y\in \partial B_x\left( 1\right) $ we have the
following estimates 
\[
A^{-1}e^{-Br\left( x\right) }< G\left( x,y\right) < Ae^{Br\left(
x\right) }.
\]
\end{lemma}

\begin{proof}
Without loss of generality we can assume that $x\in M\backslash B_p\left(
4\right) .$

Consider $\sigma $ $\left( \text{ and }\tau \right) $ the minimal geodesics
joining $p$ with $x$ (and $y$).

From now on we will use the same symbol $C$ for possibly different universal
constants (depending only on the lower bound of the Ricci curvature of $M$).

We divide the proof into two cases. In the first case, we assume that $%
d\left( x,\tau \right) >\frac 1{10},$ where $d\left( x,\tau \right)
=\min_sd\left( x,\tau \left( s\right) \right) .$ Then $G\left( x,\cdot
\right) $ defines a harmonic function on $B_z\left( 1/20\right) ,$ for any $%
z $ on $\tau ,$ and using Cheng-Yau's gradient estimate \cite{C-Y}  we conclude that there is
a constant $C$ so that
\[
\left| \nabla G\right| \left( x,z\right) \leq C G\left( x,z\right) , 
\]
for all $z$ on $\tau .$ Integrating this along $\tau $ we find that 
\[
e^{-C r\left( y\right) }G\left( x,p\right) \leq G\left( x,y\right) \leq
e^{Cr\left( y\right) }G\left( x,p\right) . 
\]
However, using the gradient estimate again we find that 
\[
e^{-Cr\left( x\right) }G\left( p,x_0\right) \leq G\left( p,x\right) \leq
e^{Cr\left( x\right) }G\left( p,x_0\right) , 
\]
for some $x_0\in \partial B_p\left( 1\right) .$ This shows that there exist
positive constants $A$ and $B$, independent of  $x$ and $y$, such
that 
\[
A^{-1}e^{-Br\left( x\right) }<G\left( x,y\right) < Ae^{Br\left(
x\right) }. 
\]
In the second case, we assume that $d\left( x,\tau \right) \leq \frac 1{10}.$
Hence there exists a point $z\in \tau $ such that $d\left( x,z\right) \leq
\frac 1{10}.$ We claim that in this case  $d\left( y,\sigma \right) >\frac 1{10}.$
Suppose the contrary, that there exists a point $w\in \sigma $ such that $%
d\left( y,w\right) \leq \frac 1{10}.$ We then have: 
\begin{eqnarray*}
d\left( p,x\right) &\leq &d\left( p,z\right) +d\left( z,x\right) \leq
d\left( p,z\right) +\frac 1{10} \\
d\left( p,y\right) &\leq &d\left( p,w\right) +d\left( w,y\right) \leq
d\left( p,w\right) +\frac 1{10}.
\end{eqnarray*}
Adding these two inequalities we get: 
\[
d\left( p,x\right) +d\left( p,y\right) \leq d\left( p,w\right) +d\left(
p,z\right) +\frac 15, 
\]
which implies that 
\[
d\left( w,x\right) +d\left( z,y\right) \leq \frac 15. 
\]
In particular, this implies that 
\[
1=d\left( x,y\right) \leq d\left( x,w\right) +d\left( w,y\right) \leq \frac
15+\frac 1{10}, 
\]
which gives a contradiction. We have thus shown that $d\left( y,\sigma
\right) >\frac 1{10}.$ Now we can conclude from the argument in the first
part of the proof that 
\[
A^{-1}e^{-Br\left( x\right) }< G\left( y,x\right) < Ae^{Br\left(
x\right) }. 
\]
The Lemma is proved.
\end{proof}

We prove a similar estimate to the Lemma above for a level set integral for $%
G$.

\begin{lemma}
\label{lem-green-int} There exist constants $C>0$ and $b>0$ such that for
any $x$ we have 
\[
\int_{l\left( s\right) }\left| \nabla G\right| \left( x,y\right) dy\leq
Ce^{br\left( x\right) }.
\]
\end{lemma}

\begin{proof}
We  claim that for $A$ and $B$ in Lemma \ref{lem-green-point} we have
\begin{equation}
L\left( 0,A^{-1}e^{-Br\left( x\right) }\right) \subset M\backslash B_x\left(
1\right) .  \label{3}
\end{equation}
Indeed, we know from Lemma \ref{lem-green-point} that 
\[
\min_{y\in \partial B_x\left( 1\right) }G\left( x,y\right) >
A^{-1}e^{-Br\left( x\right) }, 
\]
and by the maximum principle applied to $G\left( x,\cdot \right) $ on $%
B_x\left( 1\right) $ it results that 
\[
\min_{y\in B_x\left( 1\right) }G\left( x,y\right) > A^{-1}e^{-Br\left(
x\right) }. 
\]
This proves our claim that 
\[
L\left( 0,A^{-1}e^{-Br\left( x\right) }\right) \subset M\backslash B_x\left(
1\right) . 
\]
Notice that the co-area formula and property (i) for the Green's function 
 section \ref{sect-prel} imply that 
\begin{eqnarray*}
\int_{L\left( 0,A^{-1}e^{-Br\left( x\right) }\right) }\left| \nabla G\right|
^2\left( x,y\right) dy &=&\int_0^{A^{-1}e^{-Br\left( x\right) }}\left(
\int_{l\left( s\right) }\left| \nabla G\right| \left( x,y\right) dy\right) ds
\\
&=&A^{-1}e^{-Br\left( x\right) }\int_{l\left( s\right) }\left| \nabla G\right|
\left( x,y\right) dy.
\end{eqnarray*}
On the other hand, using (\ref{3}) we have that
\begin{eqnarray*}
\int_{L\left( 0,A^{-1}e^{-Br}\right) }\left| \nabla G\right| ^2\left(
x,y\right) dy &\leq &\int_{M\backslash B_x\left( 1\right) }\left| \nabla
G\right| ^2\left( x,y\right) dy \\
&=&-\int_{\partial B_x\left( 1\right) }G\left( x,y\right) \frac{\partial G}{%
\partial r}\left( x,y\right) d\sigma \left( y\right) ,
\end{eqnarray*}
where $d\sigma $ denotes the area form for $\partial B_x\left( 1\right) .$
We now use the gradient estimate and Lemma \ref{lem-green-point} to get 
\begin{eqnarray*}
\int_{L\left( 0,A^{-1}e^{-Br}\right) }\left| \nabla G\right| ^2\left(
x,y\right) dy &\leq &\int_{\partial B_x\left( 1\right) }G\left( x,y\right)
\left| \nabla G\right| \left( x,y\right) d\sigma \left( y\right)  \\
&\leq &C\int_{\partial B_x\left( 1\right) }G^2\left( x,y\right) d\sigma
\left( y\right)  \\
&\leq &Ce^{2Br\left( x\right) }.
\end{eqnarray*}
Concluding, we have that 
\[
A^{-1}e^{-Br\left( x\right) }\int_{l\left( s\right) }\left| \nabla G\left(
x,y\right) \right| dy\leq Ce^{2Br\left( x\right) },
\]
which proves the Lemma. 
\end{proof}

We now prove another result about $G\left( x,y\right) .$ We follow the proof
in \cite{Yin}.

\begin{lemma}
\label{lem-yin} For $A$ and $B$ in Lemma \ref{lem-green-point} there exists a positive constant $%
C$ such that for any $x\in M$ we have: 
\[
\int_{L\left( Ae^{Br\left( x\right) },\infty \right) }G\left( x,y\right)
dy\leq Ce^{Br\left( x\right) }.
\]
\end{lemma}

\begin{proof}

We want to prove that for any $x \in M$ 
\begin{equation}
\label{eq-inclusion}
L\left( Ae^{Br\left( x\right) },\infty \right) \subset B_x\left( 1\right) . 
\end{equation}
This is easy to see from Lemma \ref{lem-green-point}, since we know that for any $y\in \partial
B_x\left( 1\right) $ we have 
\[
G\left( x,y\right) <Ae^{Br\left( x\right) }, 
\]
and therefore by the maximum principle and the
construction of $G$ we have 
\[
\sup_{z\in M\backslash B_x\left( 1\right) }G\left( x,z\right) \leq
\sup_{y\in \partial B_x\left( 1\right) }G\left( x,y\right) <Ae^{Br\left(
x\right) }. 
\]
This proves (\ref{eq-inclusion}).

\begin{claim}
If the Ricci curvature is bounded below on a manifold $M$ with positive spectrum
 and if $G_1$ is the Dirichlet Green's function on $M$ for a ball $B_x(1)$, then 
there exists a constant $C$ such that for any $x$
\begin{equation}
\label{eq-ball}
\int_{B_x\left( 1\right) }G_1\left( x,y\right) dy<C.
\end{equation}
\end{claim}

\begin{proof}
Notice that if $H\left( z,y,t\right) $ is the Dirichlet heat kernel on $%
B_x\left( 1\right) ,$ we have 
\begin{eqnarray*}
\frac d{dt}\int_{B_x\left( 1\right) }H^2\left( x,y,t\right) dy
&=&\int_{B_x\left( 1\right) }2H\Delta _yHdy \\
&=&-2\int_{B_x\left( 1\right) }\left| \nabla _yH\right| ^2\left( x,y,t\right)
dy \\
&\leq &-2\lambda _1\left( M\right) \int_{B_x\left( 1\right) }H^2\left(
x,y,t\right) dy.
\end{eqnarray*}
This proves that 
\[
\int_{B_x\left( 1\right) }H^2\left( x,y,t\right) dy\leq e^{-2\lambda _1\left(
M\right) \left( t-1\right) }\int_{B_x\left( 1\right) }H^2\left( x,y,1\right)
dy.
\]
One can estimate the term on the right hand side as follows. First, observe that
there exists a constant $C$ independent of $x$ such that 
\[
\sup_{y\in B_x\left( 1\right) }H\left( x,y,1\right) \leq \frac C{Vol\left(
B_x\left( 1\right) \right) }.
\]
This can be proved either using the mean value inequality as in \cite{Yin}
or using the heat kernel estimates of Li and Yau \cite{L-Y} and the fact
that $H$ can be bounded from above by the heat kernel on $M$. Next we get 
\begin{eqnarray*}
\int_{B_x\left( 1\right) }H^2\left( x,y,1\right) dy &\leq &\sup_{y\in
B_x\left( 1\right) }H\left( x,y,1\right) \int_{B_x\left( 1\right) }H\left(
x,y,t\right) dy \\
&\leq &\frac C{Vol\left( B_x\left( 1\right) \right) }.
\end{eqnarray*}
We are now ready to prove our claim. Indeed, using the Cauchy-Schwarz
inequality it follows that 
\begin{eqnarray*}
\left( \int_{B_x\left( 1\right) }H\left( x,y,t\right) dy\right) ^2 &\leq
&vol\left( B_x\left( 1\right) \right) \int_{B_x\left( 1\right) }H^2\left(
x,y,t\right) dy \\
&\leq &Ce^{-2\lambda _1\left( M\right) \left( t-1\right) }.
\end{eqnarray*}
The claim now follows using that 
\[
\int_{B_x\left( 1\right) }G_1\left( x,y\right) dy=\int_0^\infty
\int_{B_x\left( 1\right) }H\left( x,y,t\right) dydt.
\]
\end{proof}

It is clear that $h=G-G_1$ is harmonic and positive on $B_x\left( 1\right) $%
. We have proved in Lemma \ref{lem-green-point} that for any $y\in \partial B_x\left( 1\right) $
we have 
\[
G\left( x,y\right) < Ae^{Br\left( x\right) }.
\]
This shows that $h< Ae^{Br\left( x\right) }$ on $\partial B_x\left(
1\right) ,$ thus by the maximum principle we know that $h<Ae^{Br\left(
x\right) }$ on $B_x\left( 1\right) .$

Then 
\begin{eqnarray}
\label{eq-bigger}
\int_{L\left( Ae^{Br\left( x\right) },\infty \right) }G\left( x,y\right) dy
&\leq &\int_{B_x\left( 1\right) }G\left( x,y\right) dy   \\
&=&\int_{B_x\left( 1\right) }G_1\left( x,y\right) dy+\int_{B_x\left(
1\right) }h\left( y\right) dy \nonumber \\
&\leq &C+vol\left( B_x\left( 1\right) \right) Ae^{Br\left( x\right) } \nonumber \\
&\leq &Ce^{Br\left( x\right) }. \nonumber
\end{eqnarray}
This finishes the proof of  Lemma \ref{lem-yin}.
\end{proof}

After this preparation, we proceed to proving Theorem \ref{thm-poisson}.

\begin{proof}[Proof of Theorem \ref{thm-poisson}]

Our goal is to show that there exist positive constants $C$ and $\alpha $
such that for any $x\in M,$ 
\begin{equation}
\left| \int_MG\left( x,y\right) f\left( y\right) dy\right| <Ce^{\alpha
r\left( x\right) }.  \label{1}
\end{equation}

Fix $x\in M$ and as above denote by $G\left( x,y\right) $ the Green's function
on $M$ with a pole at $x.$

The co-area formula (iii) and Lemma \ref{lem-green-int} yield for any $\delta >0$ and $\varepsilon >0$,
\begin{eqnarray*}
\int_{L\left( \delta \varepsilon ,\varepsilon \right) }G^{-1}\left(
x,y\right) \left| \nabla G\right| ^2\left( x,y\right) dy &=&\int_{\delta
\varepsilon }^\varepsilon \frac 1tdt\int_{l\left( t_0\right) }\left| \nabla
G\right| \left( x,y\right) dy \\
&=&\left( \int_{l\left( t_0\right) }\left| \nabla G\right| \left( x,y\right)
dy\right) \left( -\log \delta \right) \\
&\leq &Ce^{br\left( x\right) }\left( -\log \delta \right) .
\end{eqnarray*}

Notice that the integral decay of $G$ from (ii)
does not guarantee that (\ref{1}) is true, as in general it is not 
enough to balance the exponential volume growth of $M$. 
The main idea of the proof of (\ref{1}) is to use estimates on the level 
sets of $G$ instead of estimates on geodesic balls. 
We write $M$ as a disjoint union of sublevel sets of $G$, 
compute the integral of $G$  as the sum of integrals of $G$ over these sublevel sets, 
using the co-area formula and $\lambda_1\left( M\right) >0$ and estimate $f$ 
on these sublevel sets, using the given decay (\ref{eq-decay}) in order to show
 that (\ref{1}) can be bounded above by a series that converges. We first prove the following. 

\begin{claim}
There exists a constant $C$ such that for any $\delta ,\varepsilon >0$ and for any $x\in M$
\begin{equation}
\left| \int_{L\left( \delta \varepsilon ,\varepsilon \right) }G\left(
x,y\right) f\left( y\right) dy\right| \leq C\left( -\log \delta \right)
\left( \sup_{L\left( \delta \varepsilon ,\varepsilon \right) }\left| f\right|
\right) e^{br\left( x\right) }.  \label{5}
\end{equation}
\end{claim}

\begin{proof}
Let us choose the following cut off: 
\[
\phi =\chi \psi ,
\]
where we define 
\[
\chi \left( y\right) =\left\{ 
\begin{array}{c}
(\log 2)^{-1}(\log G\left( x,y\right) -\log (\frac 12\delta \varepsilon ))
\\ 
1 \\ 
(\log 2)^{-1}(\log 2\varepsilon -\log G\left( x,y\right) ) \\ 
0
\end{array}
\left. 
\begin{array}{l}
\text{on}\;\;L(\frac 12\delta \varepsilon ,\delta \varepsilon ) \\ 
\text{on}\;L\left( \delta \varepsilon ,\varepsilon \right)  \\ 
\text{on}\;\;L\left( \varepsilon ,2\varepsilon \right)  \\ 
\text{otherwise.}
\end{array}
\right. \right. .
\]
and 
\[
\psi \left( y\right) =\left\{ 
\begin{array}{c}
1 \\ 
R+1-d\left( x,y\right)  \\ 
0
\end{array}
\right. 
\begin{array}{l}
\text {on} \\ 
\text {on} \\ 
\text {on}
\end{array}
\begin{array}{l}
B_x\left( R\right)  \\ 
B_x\left( R+1\right) \backslash B_x\left( R\right)  \\ 
M\backslash B_x\left( R+1\right) 
\end{array}
.
\]
The cut off function $\chi$ will take care of the integration by parts over
 the level sets $L(\delta\epsilon,\epsilon)$ of $G$ and the cut off function $\psi$
  will take care of the integration by parts over a ball $B_x(R)$.
Then we have: 
\begin{gather*}
\lambda _1\left( M\right) \left| \int_{L\left( \delta \varepsilon
,\varepsilon \right) \cap B_x\left( R\right) }G\left( x,y\right) f\left(
y\right) dy\right| \leq \lambda _1\left( M\right) \sup_{L\left( \delta
\varepsilon ,\varepsilon \right) }\left| f\right|\times
\\
\int_{L\left(\delta \varepsilon
, \varepsilon \right) \cap B_x\left( R\right) }G\left( x,y\right) dy \\
\leq \lambda _1\left( M\right) \sup_{L\left( \delta \varepsilon ,\varepsilon
\right) }\left| f\right| \int_MG\left( x,y\right) \phi ^2\left( y\right) dy
\\
\leq \sup_{L\left( \delta \varepsilon ,\varepsilon \right) }\left| f\right|
\int_M\left| \nabla \left( G^{\frac 12}\phi \right) \right| ^2\left(
x,y\right) dy \\
\leq 2\sup_{L\left( \delta \varepsilon ,\varepsilon \right) }\left| f\right|
\left( 
\begin{array}{c}
\frac 14\int_{L\left( \frac 12\delta \varepsilon ,2\varepsilon \right)
}G^{-1}\left( x,y\right) \left| \nabla G\right| ^2\left( x,y\right) dy \\ 
+\int_MG\left( x,y\right) \left| \nabla \phi \right| ^2 (y) dy,
\end{array}
\right) 
\end{gather*}
where we have used the Poincar\'{e} inequality on $M$.
We now compute each term above. 
According to the co-area formula (iii) and Lemma \ref{lem-green-int}, we know that 
\begin{equation}
\label{eq-G-1}
\int_{L\left( \frac 12\delta \varepsilon ,2\varepsilon \right) }G^{-1}\left(
x,y\right) \left| \nabla G\right| ^2\left( x,y\right) dy\leq Ce^{br\left(
x\right) }\left( -\log \delta \right) .
\end{equation}
On the other hand, we have 
\begin{eqnarray*}
\int_MG\left( x,y\right) \left| \nabla \phi \right| ^2 &\leq &2\int_MG\left(
x,y\right) \left| \nabla \chi \right| ^2\psi ^2+2\int_MG\left( x,y\right)
\left| \nabla \psi \right| ^2\chi ^2 \\
&\leq &2\int_{L\left( \frac 12\delta \varepsilon ,2\varepsilon \right)
}G^{-1}\left( x,y\right) \left| \nabla G\right| ^2\left( x,y\right) dy \\
&&+2\int_{B_x\left( R+1\right) \backslash B_x\left( R\right) }G\left(
x,y\right) \chi ^2dy \\
&\leq &2Ce^{br\left( x\right) }\left( -\log \delta \right) +C\frac
4{\delta \varepsilon }e^{-2\sqrt{\lambda _1\left( M\right) }R}.
\end{eqnarray*}
In the last step we have used (\ref{eq-G-1}) and that since we are on the support of $\chi ,$
in particular $G\left( x,y\right) \geq \frac 12\delta \varepsilon .$ Then 
\begin{eqnarray*}
\int_{B_x\left( R+1\right) \backslash B_x\left( R\right) }G\left( x,y\right)
\chi ^2dy &\leq &\frac 2{\delta \varepsilon }\int_{B_x\left( R+1\right)
\backslash B_x\left( R\right) }G^2\left( x,y\right) dy \\
&\leq &\frac 2{\delta \varepsilon }Ce^{-2\sqrt{\lambda _1\left(
M\right) }R},
\end{eqnarray*}
using (\ref{eq-green-decay}).

We have thus proved that 
\begin{eqnarray*}
\left| \int_{L\left( \delta \varepsilon ,\varepsilon \right) \cap B_x\left(
R\right) }G\left( x,y\right) f\left( y\right) dy\right| \leq \left(
\sup_{L\left( \delta \varepsilon ,\varepsilon \right) }\left| f\right|
\right)\times\\ 
\left(
 Ce^{br\left( x\right) }\left( -\log \delta \right) +
\frac{8\bar{C}}{\delta \varepsilon }e^{-2\sqrt{\lambda _1\left( M\right) }%
R}\right) .
\end{eqnarray*}
Making now $R\rightarrow \infty $ we get that 
\[
\left| \int_{L\left( \delta \varepsilon ,\varepsilon \right) }G\left(
x,y\right) f\left( y\right) dy\right| \leq C\left( -\log \delta \right)
\left( \sup_{L\left( \delta \varepsilon ,\varepsilon \right) }\left|
f\right| \right) e^{br\left( x\right) }.
\]
This concludes the proof of the Claim. 
\end{proof}

Furthermore we have
\begin{eqnarray}
\label{eq-separate}
\left| \int_MG\left( x,y\right) f\left( y\right) dy\right| &\leq &\left|
\int_{L\left( Ae^{Br\left( x\right) },\infty \right) }G\left( x,y\right)
f\left( y\right) dy\right|  \nonumber \\
&&+\left| \int_{L\left( 0,Ae^{Br\left( x\right) }\right) }G\left( x,y\right)
f\left( y\right) dy\right| .  \label{2}
\end{eqnarray}
We will estimate each of these integrals in (\ref{eq-separate}). We start with the first integral.

Lemma \ref{lem-yin} implies that 
\begin{eqnarray}
\label{eq-1-est}
\left| \int_{L\left( Ae^{Br\left( x\right) },\infty \right) }G\left(
x,y\right) f\left( y\right) dy\right|  &\leq &\left( \sup_M\left| f\right|
\right) \int_{L\left( Ae^{Br\left( x\right) },\infty \right) }G\left(
x,y\right) dy \nonumber \\
&\leq &Ce^{Br\left( x\right) }.
\end{eqnarray}

We want to estimate the second term in equation (\ref{2}).

Recall that in Lemma \ref{lem-green-point} we proved that for any $y\in \partial B_x\left(
1\right) $ it holds 
\begin{equation}
\label{eq-below}
A^{-1}e^{-Br\left( x\right) }< G\left( x,y\right) , 
\end{equation}
where the constants $A$ and $B$ do not depend on either $x$ or $y.$ For any $z\in
M\backslash B_x\left( 1\right) $ consider $\gamma $ the minimal geodesic
joining $x$ and $z.$ Let $z_0\in \partial B_x\left( 1\right) $ be the
intersection of this geodesic with $\partial B_x\left( 1\right) .$ Then we
have, using the gradient estimate and (\ref{eq-below}) that 
\[
G\left( x,z\right) \geq G\left( x,z_0\right) e^{-Cd\left( z_0,z\right)
}\geq A^{-1}e^{-Br\left( x\right) }e^{-C_0d\left( x,z\right) }. 
\]
where $C_0>0$ is a constant not depending on $x$ or $z$. We prefer to write
this in the following equivalent way: 
\begin{equation}
G\left( x,z\right) \geq e^{-\log A-Br\left( x\right) -C_0d\left( x,z\right) },
\label{4}
\end{equation}
for all $z\in M\backslash B_x\left( 1\right) .$ The constants $A,B$ and $%
C_0$ do not depend on $x$ or $z.$ For these fixed constants let 
\[
m_0=1+\max \left\{ 2\left( B+C_0\right) r\left( x\right)+2\log A ,Br\left( x\right)
+\log A\right\} . 
\]
This choice of constant $m_0$ in particular implies that 
\[
L\left( 0,e^{-m_0}\right) \subset L\left( 0,A^{-1}e^{-Br\left( x\right) }\right)
\subset M\backslash B_x\left( 1\right) , 
\]
by (\ref{3}).

Let us write 
\begin{gather}
\left| \int_{L\left( 0,Ae^{Br\left( x\right) }\right) }G\left( x,y\right)
f\left( y\right) dy\right| \leq \left| \int_{L\left( 0,e^{-m_0}\right)
}G\left( x,y\right) f\left( y\right) dy\right|   \nonumber \\
+\left| \int_{L\left( e^{-m_0},Ae^{Br\left( x\right) }\right) }G\left(
x,y\right) f\left( y\right) dy\right|   \label{7}
\end{gather}
We now estimate $\left| \int_{L\left( 0,e^{-m_0}\right) }G\left( x,y\right)
f\left( y\right) dy\right| $ from above. First, note that from (\ref{5}) it follows that:

\begin{eqnarray}
\label{eq-G-est0}
&&\left| \int_{L\left( 0,e^{-m_0}\right) }G\left( x,y\right) f\left(
y\right) dy\right|  \nonumber \\
&=&\left| \sum_{m\geq m_0}\int_{L\left( e^{-(m+1)},e^{-m}\right) }G\left(
x,y\right) f\left( y\right) dy\right|  \nonumber \\
&\leq &\sum_{m\geq m_0}C\left( \sup_{L\left( e^{-(m+1)},e^{-m}\right)
}\left| f\right| \right) e^{br\left( x\right) }.
\end{eqnarray}

\begin{claim}
\label{claim-series}
The following series can be bounded from above by a constant $C$ independent of $x$:
\begin{equation}
\sum_{m\geq m_0}\sup_{L\left( e^{-(m+1)},e^{-m}\right) }\left| f\right|
\leq C<\infty .  \label{6}
\end{equation}
\end{claim}

\begin{proof}
On $L\left( e^{-(m+1)},e^{-m}\right) $ we know that $G\left( x,z\right)
<e^{-m}.$ The choice of constant $m_0$ implies that $L\left(
0,e^{-m_0}\right) \subset M\backslash B_x\left( 1\right) ,$ therefore using (%
\ref{4}) we get that 
\[
G\left( x,z\right) \geq e^{-\log A-Br\left( x\right) -C_0d\left( x,z\right) }, 
\]
for all $z\in L\left( 0,e^{-m_0}\right) .$ Consequently, on $L\left(
e^{-(m+1)},e^{-m}\right) $ we have 
\[
e^{-m}\geq G\left( x,z\right) \geq e^{-\log A-Br\left( x\right) -C_0d\left(
x,z\right) }. 
\]
This implies that 
\[
d\left( x,z\right) \geq \frac{m-\log A-Br\left( x\right) }{C_0}. 
\]
Moreover, we claim that 
\[
d\left( p,z\right) \geq \frac 1{2C_0}m, 
\]
for all $z\in L\left( e^{-(m+1)},e^{-m}\right) .$ To see this note that 
\begin{eqnarray*}
d\left( p,z\right) &\geq &d\left( x,z\right) -d\left( p,x\right) \\
&\geq &\frac{m-\log A-Br\left( x\right) }{C_0}-r\left( x\right) \\
&\geq &\frac m{2C_0}
\end{eqnarray*}
where the last inequality holds if and only if $m\geq 2\left( B+C_0\right)
r\left( x\right) +2\log A$, which is clearly satisfied because $m\geq m_0\geq
2\left( B+C_0\right) r\left( x\right)+2\log A$, by our choice of $m_0.$

This argument shows that 
\[
L\left( e^{-(m+1)},e^{-m}\right) \subset M\backslash B_p\left( \frac
1{2C_0}m\right) , 
\]
which we now use to estimate 
\begin{eqnarray*}
\sum_{m\geq m_0}\sup_{L\left( e^{-(m+1)},e^{-m}\right) }\left| f\right|
&\leq &\sum_{m\geq m_0}\sup_{M\backslash B_p\left( \frac 1{2C_0}m\right)
}\left| f\right| \\
&\leq &\sum_{m\geq m_0}\frac C{\left( 1+\frac 1{2C_0}m\right)
^{1+\varepsilon }}\leq C.
\end{eqnarray*}
This completes the proof of (\ref{6}).
\end{proof}

By Claim \ref{claim-series} and estimate (\ref{eq-G-est0}) we have

\begin{equation}
\label{eq-G2}
\left| \int_{L\left( 0,e^{-m_0}\right) }G\left( x,y\right) f\left( y\right)
dy\right| \leq Ce^{br\left( x\right) }. 
\end{equation}
The other term in (\ref{7}) we estimate using (\ref{5}), 
\begin{eqnarray}
\label{eq-G3}
\left| \int_{L\left( e^{-m_0},Ae^{Br\left( x\right) }\right) }G\left(
x,y\right) f\left( y\right) dy\right| &\leq &\sup_M\left| f\right|
\int_{L\left( e^{-m_0},Ae^{Br\left( x\right) }\right) }G\left( x,y\right) dy
\nonumber \\
&\leq &Ce^{br\left( x\right) }\log \left( Ae^{Br\left( x\right)
+m_0}\right) \nonumber \\
&\leq &Ce^{2br\left( x\right) }.
\end{eqnarray}
Therefore, by (\ref{7}), (\ref{eq-G2}) and (\ref{eq-G3})  we conclude that

\[
\left| \int_{L\left( 0,Ae^{Br\left( x\right) }\right) }G\left( x,y\right)
f\left( y\right) dy\right| \leq Ce^{2br\left( x\right) }.
\]
Combining  (\ref{2}), (\ref{eq-1-est})  and the previous estimate yields 
\[
\left| \int_MG\left( x,y\right) f\left( y\right) dy\right| <Ce^{\alpha
r\left( x\right) }.
\]
This concludes the proof of (\ref{1}) and hence of Theorem 1. 
\end{proof}

\section{Applications to harmonic maps and\newline
Hermitian-Einstein metrics}

\label{sect-appl}

In this section we first use the solution of the Poisson equation to
investigate the existence of harmonic maps which are homotopic to a given
map. Then we use a similar argument to study the existence of a
Hermitian-Einstein metric on a holomorphic
vector bundle over a complete manifold. These problems have been
investigated in detail in \cite{Ni, Ni1, Ni-Ren} and we rely on the
arguments from these papers. The main difference is that since we are able
to solve the Poisson equation under a lot milder assumptions, this allows us
to prove quite general existence theorems compared to the cited papers.

Results like Theorem 2 were previously proved by Li-Tam \cite{L-T} and Ding 
\cite{D}, assuming that $e\left( h\right) \in L^p\left( M\right)$ for $p$ finite. In
addition to the existence result, they were also able to prove that the
homotopy distance between $u$ and $h$ is in $L^p\left( M\right) .$ The proof
of Li and Tam was based on a heat equation approach, using ideas of Eells
and Sampson from the compact setting. Later Ni gave a different proof (he
studied the Hermitian harmonic maps, but pointed out that his proof works
for usual harmonic maps, too). It is his idea that we follow here, see also \cite{M}.

\begin{proof}[Proof of Theorem \ref{thm-harmonic}]
Let $\Omega _i$ be a compact exhaustion of $M$ with at least Lipschitz
boundary.

It is known that by a result of Hamilton there exists a sequence $u_i$ that
solves the Dirichlet homotopy problem 
\begin{eqnarray*}
\sigma \left( u_i\right)  &=&0 \\
u_i|_{\partial \Omega _i} &=&h \\
u_i &\sim &h\;\;\text{rel\thinspace \thinspace }\partial \Omega _i.
\end{eqnarray*}
In order to show convergence of  $(u_i)$ to a harmonic map 
$u:M\rightarrow N$, it is enough to prove local boundedness of the energy
density functions $e(u_i)$. To achieve this consider $\rho _i$ the homotopic
distance between $u_i$ and $h$ and $\rho _{ij}$ the homotopic distance
between $u_i$ and $u_j$. Ni has proved that a uniform bound on $\rho _{ij}$
implies a uniform bound on the energy density $e\left( u_i\right) ,$ see 
\cite{Ni} p. 344-345. His argument only uses local geometry and it is true
in our situation, too. Therefore if we can prove that the sequence $\rho
_{ij}$ is uniformly bounded on compact sets, then the convergence of $u$
follows from this argument. We now prove a uniform bound for $\rho _i$,
which by $\rho _{ij}\leq \rho _i+\rho _j$ implies a uniform bound for $\rho
_{ij}$. This is where we use the solution to the Poisson equation from
Theorem \ref{thm-poisson}. Recall that $\rho _i$ satisfy a fundamental differential
inequality
\[
\Delta \rho _i\geq -\left| \left| \sigma \left( h\right) \right| \right| ,
\]
where $\left| \left| \sigma \left( h\right) \right| \right| $ is the norm of
the tensor field. From the hypothesis we know that $\left| \left| \sigma
\left( h\right) \right| \right| $ has the right decay so we may apply
Theorem \ref{thm-poisson}. This gives a positive function $v$ that solves $\Delta v=-\left|
\left| \sigma \left( h\right) \right| \right| .$ We now use the maximum
principle to see that 
\[
\rho _i\leq v.
\]
Indeed, $v-\rho _i$ is superharmonic on $\Omega _i$ and it is positive on $%
\partial \Omega _i.$ This proves that $\rho _i$ is uniformly bounded on a
fixed compact set $K\subset M,$ therefore the energy density $e\left(
u_i\right) $ is uniformly bounded on $K.$ This proves the existence of $u$.
Since by Theorem \ref{thm-poisson} the function $v$ has at most the exponential growth,
it follows that the homotopic distance between $u$ and $h$ grows at most exponentially.
\end{proof}

We now discuss another application of Theorem \ref{thm-poisson}, which
concerns the existence of Hermitian-Einstein metrics on holomorphic vector
bundles over a complete K\"{a}hler manifold, stated in Theorem \ref{thm-he}.

The existence of a Hermitian Einstein metric on a compact manifold is
related to stability of the vector bundle, as it is well known from the work
of Donaldson and Uhlenbeck-Yau. On complete manifolds Ni and Ni-Ren \cite
{Ni1, Ni-Ren} have showed that if there exists a Hermitian metric $H_0$ such
that $\left| \left| \Lambda F_{H_0}-\lambda I\right| \right| \in L^p\left(
M\right) $ for some $p\geq 1$ finite and $M$ has positive spectrum then there
exists a metric $H$ such that $\Lambda F_H=\lambda I.$ As in the case of
harmonic maps, we can modify their argument to obtain Theorem \ref{thm-he}.

\begin{proof}[Proof of Theorem \ref{thm-he}]
The argument here is similar to the argument for harmonic maps. Here
we follow \cite{Ni1}, and use again our solution from Theorem \ref{thm-poisson} and the
maximum principle to prove certain $C^0$ estimates. We
discuss here the special case $\lambda =0,$ however the general case follows the
same. By a result of Donaldson there exists a hermitian metric $H_i$ on $%
\Omega _i$ such that 
\begin{eqnarray*}
\Lambda F_{H_i} &=&0\;\;\text{on\thinspace \thinspace }\Omega _i \\
H_i|_{\partial \Omega _i} &=&H_0.\;
\end{eqnarray*}

The goal is to prove that we can pass to a limit and obtain a solution on $M$
and for this we establish a priori estimates. Recall the following distance
functions introduced by Donaldson
\begin{eqnarray*}
\tau _i\left( x\right)  &=&\tau \left( H_i,H_0\right) =tr\left(
H_iH_0^{-1}\right)  \\
\sigma _i\left( x\right)  &=&\sigma \left( H_i,H_0\right) =tr\left(
H_iH_0^{-1}\right) +tr\left( H_0H_i^{-1}\right) -2rank\left( E\right) .
\end{eqnarray*}
In order to prove that $H_i$ has a subsequence which converges uniformly on
compact subsets of $M$ we need to establish $C^0$ and $C^1$ estimates. The $%
C^0$ estimates is based on the solution to the Poisson equation found in
Theorem 1. Recall that we have a Bochner type inequality for the distance $%
\tau _i$, see \cite{Siu, Ni1}. Let $f_i=\log \tau _i-\log k,$ where $k$ is
the rank of $E.$ Then
\begin{eqnarray*}
\Delta f_i &\geq &-\left| \left| \Lambda F_{H_0}\right| \right|
\;\;\text{on}\;\;\Omega _i \\
f_i &=&0\;\;\text{on}\;\;\partial \Omega _i.
\end{eqnarray*}
Using the maximum principle it is easy to see that 
\[
f_i\leq v,
\]
where $v$ is the solution of 
\[
\Delta v=-\left| \left| \Lambda F_{H_0}\right| \right| .
\]

Such a solution exists and has at most exponential growth since the norm of $%
\Lambda F_{H_0}$ has the decay as in Theorem \ref{thm-poisson}. This shows that $\tau _i\leq
ke^v$ and therefore we get a bound for $\sigma _i$%
\[
\sigma _i\leq 2ke^v-2k.
\]
This establishes the $C^0$ estimates and now the $C^1$ estimates follow the
same as in \cite{Ni1}, p.690-691. Since the argument for these $C^1$
estimates only uses local geometry it is valid in our situation as well.
This concludes the proof.
\end{proof}

\subsection{The Ricci flow}

In order to Prove Theorem \ref{thm-rf1} we want to find a bounded solution
to the Poisson equation 
\[
\Delta u=(R+1), 
\]
at time $t=0$, where $\Delta $ is the Laplacian taken with respect to metric 
$g_0$.

The following result holds in arbitrary dimension.

\begin{proposition}
\label{prop-thm1} Let $(M,g_0)$ be a Cartan Hadamard manifold of dimension $n.$
Assume that its Ricci curvature satisfies 
\[
-a^2\le \mathrm{Ric}_M\le -b^2<0. 
\]
Then for any $f\in C^\infty (M)$ having a decay 
\[
|f(x)|\le \frac C{(1+r(x))^{1+\varepsilon }} 
\]
where $\varepsilon >0$ and $r(x)$ is a distance from $x$ to a fixed point $p$%
, there exists a unique solution to $\Delta u=f$ having decay 
\[
|u(x)|\le \frac{\tilde{C}}{(1+r(x))^\varepsilon }. 
\]
\end{proposition}

\begin{proof}
Let us observe first that $M$ has positive spectrum. Since it is simply 
connected and its Ricci curvature is bounded above by -$b^2$, 
the Laplacian comparison theorem (see \cite {X}, Theorem 2.15) implies that 
\[
\Delta r\geq b\coth br.
\]
Thus, if $f$ is any compactly supported smooth function on $M$ it results that 
\begin{eqnarray*}
b\int_Mf^2\left( x\right) \coth \left( br\left( x\right) \right)  &\leq
&\int_Mf^2\left( x\right) \Delta r=-2\int_Mf\nabla f\cdot \nabla r \\
&\leq &2\int_M\left| f\right| \left| \nabla f\right|  \\
&\leq &\frac b2\int_Mf^2+\frac 2b\int_M\left| \nabla f\right| ^2.
\end{eqnarray*}
This implies that 
\[
\frac{b^2}4\int_Mf^2\leq \int_M\left| \nabla f\right| ^2,
\]
and from here we infer that 
\[
\lambda _1\left( M\right) \geq \frac{b^2}4>0.
\]
Consider the sequence $u_i : B_p(R_i) \to \mathbb{R}$, for $R_i \to \infty$,
such that
$$u_i(x) = -\int_{B_p(R_i)} G_i(x,y)f(y)\, dy,$$
where $G_i(x,y)$ is the Dirichlet Green's function on $B_p(R_i)$. Then
\begin{eqnarray*}
\Delta u_i &=& f, \,\,\, \mbox{on} \,\,\, B_p(R_i) \\
u_i &=& 0, \,\,\, \mbox{on} \,\,\, \partial B_p(R_i).
\end{eqnarray*}
Our goal is to show that we have uniform estimates on $u_i$ so that
passing to a limit as $i\to\infty$ yields to a global solution
to the Poisson equation $\Delta u = f$ with a desired decay at infinity.
Observe that by Laplacian comparison theorem
$$\Delta r \ge b\coth(br) \ge b,$$
where $r(\cdot) = d(\cdot,p)$. We have
\begin{eqnarray*}
\Delta \frac{1}{r^{\varepsilon}} &=& -\frac{\varepsilon}{r^{1+\varepsilon}}\Delta r + 
\frac{\varepsilon(\varepsilon + 1)}{r^{2+\varepsilon}} \\
&\le& -\frac{b\varepsilon}{r^{1+\varepsilon}} +   
\frac{\varepsilon(\varepsilon + 1)}{r^{2+\varepsilon}}.
\end{eqnarray*}
By choosing $r \ge r_0 \ge 1$ with $r_0$ sufficiently big e.g. $r_0=max(\frac{2(1+\varepsilon)}{b},1)$
we get
\begin{equation}
\label{eq-alpha}
\Delta r^{-\varepsilon} \le -\frac{\alpha}{r^{1+\varepsilon}},
\end{equation}
where $\alpha= \frac{\varepsilon b}{2}>0$. We can use a family of balls
$\{B_p(R_i)\}$ as an exhaustion family of $M$ by compact sets
in the construction of Green's function $G(x,y)$. Therefore,
by the maximum principle we have
$$G_i(x,y) \le G(x,y), \,\,\, \mbox{for every} \,\,\, x,y\in B_p(R_i).$$
We have the following,
\begin{eqnarray*}
|u_i(x)| &\le& \int_{B_p(R_i)} G_i(x,y)|f(y)|\, dy \\
&\le& \int_{B_p(R_i)} G(x,y)|f(y)|\, dy \\
&\le& \int_M G(x,y)|f(y)|\, dy < \infty,
\end{eqnarray*}
where the finiteness of the last term follows by Theorem
\ref{thm-poisson}. This shows $u_i$ is uniformly bounded on compact
sets of $M$ and by standard elliptic estimates and Arzela-Ascoli
theorem we can extract a subsequence of $u_i$ converging uniformly on
compact sets to a smooth limit $u(x)$ that satisfies $\Delta u = f$ on
$M$. Without loss of generality we still denote this convergent subsequence with
$u_i$.
Moreover, to establish a better decay of $u(x)$ at infinity, let $A$ be such that

$$A > r_0^{\varepsilon} \cdot \sup_i\sup_{x\in \partial B_p(r_0)}|u_i(x)|$$
and $A >
\frac{C}{\alpha}$, where $C$ is the constant from (\ref{eq-decay}) and
$\alpha$ is the same as above. Consider 
$$v_i = \frac{A}{r^{\varepsilon }} - u_i.$$ 
Then by (\ref{eq-alpha}) and $r_0 > 1$,
\begin{eqnarray}
\label{eq-lap}
\Delta v_i &\le& - A\frac{\alpha}{r^{\varepsilon + 1}}
- f \nonumber\\
&\le& -\frac{A\alpha}{r^{\varepsilon+1}} + \frac{C}{r^{\varepsilon + 1}} \le 0.
\end{eqnarray}
We also have
$$v_i|_{\partial B_p(r_0)} = Ar_0^{-\varepsilon} - u_i|_{\partial B_p(r_0)} > 0,$$
$$v_i|_{\partial B_p(R_i)} > 0 \,\,\, \mbox{since} \,\,\, u_i|_{\partial B_p(R_i)} = 0.$$
By the maximum principle applied to $v_i$ satisfying (\ref{eq-lap}), we get
$v_i \ge 0$ on $B_p(R_i)\backslash B_p(r_0)$ and therefore
$$u_i \le \frac{A}{r^{\varepsilon}}, \,\,\, \mbox{on}
\,\,\, B_p(R_i)\backslash B_p(r_0).$$
Applying the same arguments to $-u_i$ we get
$$|u_i| \le \frac{ A}{r^{\varepsilon}}, \,\,\, \mbox{on}
\,\,\, B_p(R_i)\backslash B_p(r_0).$$
Letting $i\to \infty$ we conclude the proof.
\end{proof}

\begin{proof}[Proof of Theorem \ref{thm-rf1}]
By Proposition \ref{prop-thm1} we can find a bounded potential
solving
$$\Delta u = R+1.$$
By Theorem $1.1$ in \cite{Ch} it immediatelly follows the Ricci flow
has a long time existence on $M$ and it converges, as $t\to \infty$,
uniformly on compact subsets, to a K\"ahler Einstein metric with 
constant negative curvature. 
\end{proof}

\end{document}